\documentclass[reqno,14pt,centertags]{amsart}
\usepackage{amsmath,amsthm,amscd,amssymb,latexsym,upref,stmaryrd}

\usepackage{hyperref}
\newcommand*{\mailto}[1]{\href{mailto:#1}{\nolinkurl{#1}}}

\newcommand{\B}{$\hfill\Box$}
\newcommand{\al}{\alpha}

\newcommand{\de}{\delta}
\newcommand{\la}{\lambda}
\newcommand{\om}{\omega}

\newcommand{\vv}{\varphi}
\newcommand{\iy}{\infty}

\newtheorem{theorem}{Theorem}[section]
\newtheorem{lemma}[theorem]{Lemma}

\theoremstyle{definition}

\begin{document}
\title[Recovering Dirac Operator with Nonlocal Boundary Conditions]
{Recovering Dirac Operator with Nonlocal Boundary Conditions}

\author[C.~F.~Yang]{Chuan-Fu Yang}
\address{Department of Applied Mathematics, Nanjing University of
Science and Technology, Nanjing, 210094, Jiangsu, People's
Republic of China}
\email{\mailto{chuanfuyang@mail.njust.edu.cn}}

\author[V.~Yurko]{Vjacheslav Yurko}
\address{Department of Mathematics, Saratov University,
Astrakhanskaya 83, Saratov 410012, Russia}
\email{\mailto{yurkova@info.sgu.ru}}

\date{\today}

\begin{abstract}
In this paper inverse problems for Dirac operator with nonlocal conditions are considered.
Uniqueness theorems of inverse problems from the Weyl-type function and spectra
are provided, which are generalizations of the well-known Weyl function and Borg's inverse
problem for the classical Dirac operator.
\end{abstract}

\maketitle

\section{Introduction}

Problems with nonlocal conditions arise in various fields
of mathematical physics \cite{d,g,i,1,yin}, biology and biotechnology \cite{n,s}, and
in other fields. Nonlocal conditions come up when
value of the function on the boundary is connected to values inside the domain.
Recently problems with nonlocal conditions are paid much attention for them
in the literature.

In this paper we study inverse spectral problems for Dirac operator
\begin{align}\label{1}
By'+\Omega(x)y=\lambda y,\quad x\in (0,T),
\end{align}
and with nonlocal linear conditions
\begin{align}\label{2}
U_j(y):=\int_0^T y(t)^td\sigma_j(t)=0,\quad j=1,2.
\end{align}
Here
\begin{align*}
B=\left(\begin{array}{cc}
0&1\\
-1&0\end{array}\right),\ \ \Omega(x)=\left(\begin{array}{cc}
p(x)&q(x)\\
q(x)&-p(x)\end{array}\right),\ \ y(x)=\left(\begin{array}{c}
y_1(x)\\
y_2(x)\end{array}\right),
\end{align*}
functions $p(x)$ and $q(x)$ are complex-valued and absolutely continuous functions in $(0, T)$, and $\lambda$ is a spectral
parameter, vector-valued functions $\sigma_j(t)=\left(\begin{array}{c}
\sigma_{j,1}(t)\\
\sigma_{j,2}(t)\end{array}\right)$ are
complex-valued functions of bounded variations and are continuous from
the right for $t>0.$ There exist finite limits $h_{j,i}:=\sigma_{j,i}(+0)-
\sigma_{j,i}(0),\ i=1,2.$ Linear forms $U_j$ in (\ref{2}) can be written as forms
\begin{align}\label{3}
U_j(y):=h_{j,1}y_1(0)+h_{j,2}y_2(0)+\int_0^T y(t)^td\sigma_{j0}(t),\quad j=1,2,
\end{align}
where vector-valued functions $\sigma_{j0}(t)$ in (\ref{3}) are complex-valued functions of bounded
variations and are continuous from the right for $t\ge 0$, and $|h_{1,1}|+|h_{1,2}|\ne 0.$

A complex number $\lambda_0$ is called an eigenvalue of the problem (\ref{1}) and (\ref{2}) if
equation (\ref{1}) with $\lambda=\lambda_0$ has a nontrivial solution $y_0(x)$
satisfying conditions (\ref{2}); then $y_0(x)$ is called the eigenfunction of the problem
(\ref{1}) and (\ref{2}) corresponding to the eigenvalue $\lambda_0$. The number of linearly independent solutions of the
problem (\ref{1}) and (\ref{2}) for a given eigenvalue $\lambda_0$ is called the multiplicity of $\lambda_0$.

Classical inverse problems for
Eq.(\ref{1}) with two-point boundary conditions have been studied
fairly completely in many works (see \cite{GAS,GDZ,GUS,LSA} and the
references therein). The theory of nonlocal inverse spectral problems
now is only at the beginning because of its complexity. Results
of the inverse problem for various nonlocal operators can be
found in \cite{7,8,9,10,11,12,13,14,15}.

In this work by using Yurko's ideas of the method of spectral mappings \cite{5} we prove uniqueness theorems for the
solution of the inverse spectral problems for Eq.(\ref{1}) with nonlocal
conditions (\ref{2}).

\section{main results}

Let $X_k(x,\la)$ and $Z_k(x,\la),$ $k=1,2,$ be the solutions of Eq.(\ref{1})
with the initial conditions
$$
X_1(0,\la)=\left(\begin{array}{c}
1\\
0\end{array}\right)=Z_1(T,\la),\ X_2(0,\la)=\left(\begin{array}{c}
0\\
1\end{array}\right)=Z_2(T,\la).
$$
Denote by $L_0$ the boundary value problem (BVP) for Eq.(\ref{1}) with the conditions
$$
U_1(y)=U_2(y)=0,
$$
and $\om(\la):=\det[U_j(X_k)]_{j,k=1,2}$, and assume that $\om(\la)\not\equiv 0.$
The function $\om(\la)$ is an entire function of exponential type with order $1,$ and its zeros
$\Xi:=\{\xi_n\}_{n\in \mathbb{Z}}$ (counting multiplicities) coincide with the eigenvalues of $L_0$. The function $\om(\la)$
is called the characteristic function for $L_0$.

Denote $V_j(y):=y_j(T),\; j=1,2.$ Consider the BVP $L_j$, $j=1,2,$ for
Eq.(1) with the conditions $U_j(y)=V_1(y)=0.$ The eigenvalue sets $\Lambda_j:=
\{\la_{nj}\}_{n\in \mathbb{Z}}$ (counting multiplicities) of the BVP $L_j$ coincide with the zeros of the
characteristic function $\Delta_j(\la):=\det[U_j(X_k), V_1(X_k)]_{k=1,2}$.

For $\lambda\neq \lambda_{n1}$, let $\Phi(x,\la):=\left(\begin{array}{c}
\Phi_1(x,\la)\\
\Phi_2(x,\la)\end{array}\right)
$ be the solution of Eq.(\ref{1}) under the conditions $U_1(\Phi)=1$ and
$V_1(\Phi):=\Phi_1(T,\la)=0.$ Denote Weyl-type function $M(\la):=U_2(\Phi).$ It is known \cite{LSA} that for Dirac operator with
classical two-point separated boundary conditions, the specification of the
Weyl function uniquely determines the function $\Omega(x)$. However, in the case with
nonlocal conditions, it is not true; the specification of the Weyl-type
function $M(\la)$ does not uniquely determine the function $\Omega(x)$ (see counterexamples
in Section 5). For the nonlocal conditions the inverse problem is formulated as follows.

Throughout this paper the functions $\sigma_{ji}(t)$ are known a priori.
And condition $S$: $\Lambda_1\cap\Xi=\emptyset$.

Let us formulate a uniqueness theorem. For this purpose,
together with $\Omega(x)$ we consider another $\tilde{\Omega}(x)$, and we agree that
if a certain symbol $\al$ denotes an object related to $\Omega(x)$, then $\tilde\al$
will denote an analogous object related to $\tilde{\Omega}(x)$.

{\bf Theorem 1. }{\it Let condition $S$ be true. If $M(\la)=\tilde M(\la)$
and $\om(\la)=\tilde\om(\la),$ then $\Omega(x)=\tilde \Omega(x)$ on $(0,T).$}

Consider the BVP $L_{11}$ for Eq.(\ref{1}) with the conditions $U_1(y)=0=V_2(y):=y_2(T).$
The eigenvalue set $\Lambda_{11}:=\{\la_{n1}^1\}_{n\in \mathbb{Z}}$ of the BVP $L_{11}$
coincide with the zeros of the characteristic function $\Delta_{11}(\la):=
\det[U_1(X_k), V_2(X_k)]_{k=1,2}$. Clearly, $\{\lambda_{n1}\}_{n\in \mathbb{Z}}\cap\{\lambda_{n1}^1\}_{n\in \mathbb{Z}}=\emptyset.$

{\bf Theorem 2. }{\it If $\la_{n1}=\tilde\la_{n1}, \la_{n1}^1=\tilde\la_{n1}^1$,
$n\in \mathbb{Z},$ then $\Omega(x)=\tilde \Omega(x)$ on $(0,T).$}

This theorem is an analogue of the well-known Borg's inverse
problem \cite{16} for Sturm-Liouville operators with classical two-point
separated boundary conditions.

\section{Lemmas}

Denote $\Lambda^{\pm}:=\{\lambda:\pm\mbox{Im}\lambda\geq 0\}$.

\begin{lemma}(See \cite{LSA}.)\label{yl1}
For $\lambda\in \Lambda^{\pm}$, $|\lambda|\rightarrow\infty$, the following asymptotic formulas hold:
$$
X_1(x,\lambda)=\left(\begin{array}{c}
\cos(\lambda x)\\
\sin(\lambda x)\end{array}\right)\left[1+O\left(\frac{1}{\lambda}\right)\right],
$$
$$
X_2(x,\lambda)=\left(\begin{array}{c}
-\sin(\lambda x)\\
\cos(\lambda x)\end{array}\right)\left[1+O\left(\frac{1}{\lambda}\right)\right].
$$
\end{lemma}

By Lemma \ref{yl1} there exists a fundamental system of
solutions $Y_k(x,\lambda)\ (k=1,2)$ of Eq.(1) such that for $\lambda\in\Lambda^{\pm}$, $|\lambda|\to\iy$:
\begin{equation}\label{4}
\begin{array}{l}
Y_1(x,\lambda)\!=\!\exp(i\lambda x)\left[\left(\begin{array}{c}
1\\
-i\end{array}\right)+O\left(\frac{1}{\lambda}\right)\right],
\\
Y_2(x,\lambda)\!=\!\exp(-i\lambda x)\left[\left(\begin{array}{c}
1\\
i\end{array}\right)+O\left(\frac{1}{\lambda}\right)\right].
\end{array}
\end{equation}
and
\begin{align}\label{5}
\det[Y_1(x,\lambda),Y_2(x,\lambda)]=2i\left[1+O\left(\frac{1}{\lambda}\right)\right].
\end{align}

{\bf Lemma 1. }(See \cite{yy}) {\it Let $\{W_k(x,\la)\}_{k=1,2}$ be a fundamental system of
solutions of Eq.(\ref{1}), and let $Q_j(y),\; j=1,2,$ be linear forms. Then}
\begin{align}\label{6}
\det[Q_j(W_k)]_{k,j=1,2}=\det[Q_j(X_k)]_{k,j=1,2}
\det[W_k^{(\nu-1)}(x,\la)]_{k,\nu=1,2}\,.
\end{align}

It follows from (\ref{5})-(\ref{6}) that
\begin{align}\label{7}
\det[Q_j(Z_k)]_{k,j=1,2}=\det[Q_j(X_k)]_{k,j=1,2},
\end{align}
and
\begin{align}\label{8}
\det[Q_j(Y_k)]_{k,j=1,2}=2i(1+O(\lambda^{-1})) \det[Q_j(X_k)]_{k,j=1,2}\,.
\end{align}

Introduce the functions
$$
\vv(x,\la)=U_1(X_1)X_2(x,\la)-U_1(X_2)X_1(x,\la),
$$
$$
\theta(x,\la)=U_2(X_2)X_1(x,\la)-U_2(X_1)X_2(x,\la),
$$
$$
\psi(x,\la)=V_1(X_2)X_1(x,\la)-V_1(X_1)X_2(x,\la).
$$
Then
$$
U_1(\vv)=0,\;U_2(\vv)=\om(\la),\;V_1(\vv)=\Delta_1(\la),\;
V_2(\vv)=\Delta_{11}(\la),
$$
$$
U_1(\theta)=\om(\la),\; U_2(\theta)=0,\; V_1(\theta)=-\Delta_2(\la),
$$
$$
U_j(\psi)=\Delta_j(\la),\; V_1(\psi)=0,\; V_2(\psi)=-1.
$$
Moreover, by calculation, Eqs.(\ref{6})-(\ref{7}) yield that
\begin{equation}\label{9}
\begin{array}{l}
\det[\theta(x,\la), \vv(x,\la)]\!=\!\om(\la),\
\det[\psi(x,\la), \vv(x,\la)]\!=\!\Delta_1(\la),
\end{array}
\end{equation}
\begin{align}\label{10}
\Delta_1(\la)=-U_1(Z_2),\;\Delta_2(\la)=-U_2(Z_2),\;\Delta_{11}(\la)=U_1(Z_1).
\end{align}
Note that the functions $\Phi, \psi, \vv$ and $\theta$ are all
the solutions of Eq.(\ref{1}) with some conditions, comparing boundary conditions on $\Phi, \psi, \vv$ and $\theta,$ we arrive at
\begin{align}\label{11}
\Phi(x,\la)=\frac{\psi(x,\la)}{\Delta_1(\la)},
\end{align}
\begin{align}\label{12}
\Phi(x,\la)=\frac{1}{\om(\la)}\,\Big(\theta(x,\la)+\frac{\Delta_2(\la)}{\Delta_1(\la)}\vv(x,\la)\Big).
\end{align}
Hence,
\begin{align}\label{13}
M(\la):=U_2(\Phi)=\frac{\Delta_2(\la)}{\Delta_1(\la)},
\end{align}
\begin{align}\label{14}
\det[\Phi(x,\la), \vv(x,\la)]=1.
\end{align}

Let $v_1(x,\la)=\left(\begin{array}{c}
v_{11}(x,\la)\\
v_{12}(x,\la)\end{array}\right)$ and $v_2(x,\la)=\left(\begin{array}{c}
v_{21}(x,\la)\\
v_{22}(x,\la)\end{array}\right)$ be the solutions of Eq.(\ref{1}) with
the conditions
$$
v_1(T,\la)=\left(\begin{array}{c}
1\\
0\end{array}\right),
$$
and
$$
v_{22}(T,\la)=1,\ U_1(v_2)=0.
$$
Obviously,
\begin{equation}\label{15}
\begin{array}{l}
v_1(x,\la)=Z_1(x,\la),\; v_2(x,\la)=Z_2(x,\la)+N(\la)Z_1(x,\la),\\
\det[v_1(x,\la),v_2(x,\la)]=1,
\end{array}
\end{equation}
where
\begin{align}\label{16}
N(\la)=\frac{\Delta_1(\la)}{\Delta_{11}(\la)}=-\frac{U_1(Z_2)}{U_1(Z_1)}\,.
\end{align}

Denote
$$
U_1^a(y):=\int_0^a y(t)d\sigma_1(t),\quad a\in(0,T].
$$
Clearly, $U_1=U_1^T,$ and if $\sigma_1(t)\equiv C$ (constant) for $t\ge a,$
then $U_1=U_1^a.$

For sufficiently small $\de>0,$ we denote
$$
\Pi_\de:=\{\lambda:\; \mbox{arg}\,\lambda\in[\de,\pi-\de]\},\ \
G_\de:=\{\lambda:\; |\lambda-\la_{n1}|\ge\de,\;\;\forall n\in \mathbb{Z}\},
$$
and
$$
G_\de':=\{\lambda:\; |\lambda-\la_{n1}^1|\ge\de,\;\;\forall n\in \mathbb{Z}\},
$$
where $\la_{n1}\in \Lambda_1$ and $\la_{n1}^1\in \Lambda_{11}$.

{\bf Lemma 2. }{\it For $\lambda\in\Pi_\de,\; |\lambda|\to\iy,$ there hold
\begin{align}\label{17}
\Phi(x,\la)=
(h_{1,1}-ih_{1,2})^{-1}\exp(i\lambda x)\left[\left(\begin{array}{c}
1\\
-i\end{array}\right)+o(1)\right],\;x\in[0,T),
\end{align}
\begin{align}\label{170}
\Psi(x,\la)=\frac{1}{2i}\exp(-i\lambda (T-x))\left[\left(\begin{array}{c}
1\\
-i\end{array}\right)+o(1)\right],\;x\in[0,T),
\end{align}
\begin{align}\label{18}
v_1(x,\la)=\frac{1}{2}\exp(-i\lambda(T-x))\left[\left(\begin{array}{c}
1\\
-i\end{array}\right)+O\left(\frac{1}{\lambda}\right)\right],\; x\in[0,T),
\end{align}
\begin{equation}\label{19}
\begin{array}{l}
\Delta_1(\la)=\frac{h_{1,1}-ih_{1,2}}{2i}\,\exp\left(-i\lambda T\right)[1+o(1)],\\
\Delta_{11}(\la)=\frac{h_{1,1}-ih_{1,2}}{2}\,\exp\left(-i\lambda T\right)[1+o(1)].
\end{array}
\end{equation}

Let $\sigma_1(t)\equiv C$ (constant) for $t\ge a$ (i.e. $U_1=U_1^a$). Then for
$\lambda\in\Pi_\de,\; |\lambda|\to\iy,$} we get
\begin{align}\label{20}
&\vv(x,\la)=\frac{h_{1,1}-ih_{1,2}}{2i}\exp(-i\lambda x)\left[\left(\begin{array}{c}
1\\
i\end{array}\right)+o(1)+O(\exp(i\lambda(2x-a)))\right],\; x\in(0,T],
\end{align}
\begin{align}\label{21}
&v_2(x,\la)=\exp(i\lambda(T-x))\left[\left(\begin{array}{c}
-i\\
1\end{array}\right)+o(1)+O(\exp(i\lambda(2x-a)))\right]\; x\in[0,T).
\end{align}

{\it Proof.} The function $\Phi(x,\la)$ can be expressed as
\begin{align}\label{22}
\Phi(x,\la)=A_1(\la)Y_1(x,\lambda)+A_2(\la)Y_2(x,\lambda),
\end{align}
together with $U_1(\Phi)=1$ and $V_1(\Phi):=\Phi_1(T,\lambda)=0,$ which yields that
\begin{align}\label{23}
A_1(\la)U_1(Y_1)+A_2(\la)U_1(Y_2)=1,\; A_1(\la)V_1(Y_1)+A_2(\la)V_1(Y_2)=0.
\end{align}
Using (4), one gets that for $\lambda\in\Pi_\de,\; |\lambda|\to\iy$:
\begin{align}\label{24}
U_1(Y_1)=(h_{1,1}-ih_{1,2})(1+o(1)),\; U_1(Y_2)=O(\exp(-i\lambda T)),
\end{align}
\begin{equation}\label{25}
\begin{array}{l}
V_1(Y_1)=\exp\left(i\lambda T\right)(1+O(\lambda^{-1})),\
V_1(Y_2)=\exp\left(-i\lambda T\right)(1+O(\lambda^{-1})).
\end{array}
\end{equation}
Solving linear algebraic system (\ref{23}) by using (\ref{24})-(\ref{25}), we obtain
$$
A_1(\lambda)=(h_{1,1}-ih_{1,2})^{-1}(1+o(1)),\; A_2(\lambda)=O(\exp(2i\lambda T)).
$$
Substituting these relations into (\ref{22}), we have proved (\ref{17}).
Formulas (\ref{18})-(\ref{21}) can be proved similarly, and are omitted. \B

By the well-known method (see, for example, \cite{4}) the
following estimates hold for $x\in(0,T),\; \lambda\in \Lambda^+:$
\begin{align}\label{26}
v_1(x,\la)=O\left(\exp\left(-i\lambda(T-x)\right)\right),
\end{align}
\begin{align}\label{27}
\Phi(x,\la)=O\left(\exp\left(i\lambda x\right)\right),\quad \rho\in G_\de.
\end{align}
Moreover, if $\sigma_1(t)\equiv C$ (constant) for $t\ge a$ (i.e. $U_1=U_1^a$),
then for $x\ge a/2,\; \lambda\in \Lambda^+$:
\begin{align}\label{28}
\vv(x,\la)=O\left(\exp\left(-i\lambda x\right)\right),
\end{align}
\begin{align}\label{29}
v_2(x,\la)=O\left(\exp\left(i\lambda(T-x)\right)\right),\quad \rho\in G'_\de.
\end{align}

\section{Proofs of Theorems}

{\textbf{Proof of Theorem 2}}
We know that the characteristic function $\Delta_1(\la)$ of the BVP $L_1$ is an entire function of order one with respect to $\lambda$.
Following the theory of Hadamard's factorization (see \cite{AHL}), $\Delta_1(\lambda)$ can be expressed as an
infinite product as
$$
\Delta_1(\lambda)=c_1e^{a_1\lambda}\prod_{n\in \mathbb{Z}}\left(1-\frac{\lambda}{\lambda_{n1}}\right)
e^{\frac{\lambda}{\lambda_{n1}}+\frac{1}{2}(\frac{\lambda}{\lambda_{n1}})^2+\cdots+\frac{1}{p}(\frac{\lambda}{\lambda_{n1}})^p},
$$
where $\lambda_{n1}$ are the eigenvalues of the problem $L_1$, $p$ is the genus of $\Delta_1(\lambda)$,
$c_1$ and $a_1$ are constants. Since for
the order $\rho$ of $\Delta_1(\lambda)$, $p\leq \rho\leq p+1$ (see \cite{AHL}), and $\Delta_1(\lambda)$ is
an entire function of exponential type with order $1,$ we find that the genus of $\Delta_1(\lambda)$ is $0$ or $1$ (that is, $p=0\vee 1$). Thus
$\Delta_1(\lambda)$ can be rewritten by
$$
\Delta_1(\lambda)=c_1e^{a_1\lambda}\prod_{n\in \mathbb{Z}}\left(1-\frac{\lambda}{\lambda_{n1}}\right)
e^{\frac{\lambda}{\lambda_{n1}}p}.
$$

Since $\Delta_1(\lambda)$ and $\tilde{\Delta}_1(\lambda)$ are both entire functions of order one with respect
to $\lambda$, and $\lambda_{n1}=\tilde{\lambda}_{n1}$ for all $n\in \mathbb{Z}$, by the Hadamard's factorization theorem, we may suppose
(the case when $\Delta_1(0)=0$ requires minor modifications)
$$
\Delta_1(\lambda)=c_1e^{a_1 \lambda}\prod_{n\in \mathbb{Z}}\left(1-\frac{\lambda}{\lambda_{n1}}\right)e^{\frac{\lambda}{\lambda_{n1}}p}
$$
and
$$
\tilde{\Delta}_1(\lambda)=\tilde{c}_1e^{\tilde{a}_1 \lambda}\prod_{n\in \mathbb{Z}}\left(1-\frac{\lambda}{\lambda_{n1}}\right)e^{\frac{\lambda}{\lambda_{n1}}\tilde{p}},
$$
for some constants $c_1,\tilde{c}_1,a_1,\tilde{a}_1$ and $p,\tilde{p}$, which can be determined from
the asymptotics.
From this we get for all $\lambda\in \mathbb{C}$
$$
\frac{\tilde{\Delta}_1(\lambda)}{\Delta_1(\lambda)}=\frac{\tilde{c}_1}{c_1}
e^{\left[(\tilde{a}_1-a_1)+(\tilde{p}-p)\sum_{n\in \mathbb{Z}}\frac{1}{\lambda_{n1}}\right]\lambda}.
$$
The expression (\ref{19}) implies that
$$
\frac{\tilde{\Delta}_1(\la}{\Delta_1(\la)}=1+o(1)\equiv\frac{\tilde{c}_1}{c_1}
e^{\left[(\tilde{a}_1-a_1)+(\tilde{p}-p)\sum_{n\in \mathbb{Z}}\frac{1}{\lambda_{n1}}\right]\lambda},
$$
which yields that
$$
(\tilde{a}_1-a_1)+(\tilde{p}-p)\sum_{n\in \mathbb{Z}}\frac{1}{\lambda_{n1}}=0,\ \tilde{c}_1=c_1.
$$
Consequently,
$\Delta_1(\la)\equiv\tilde\Delta_1(\la).$ Analogously, from $\lambda_{n1}^1=\tilde{\lambda}_{n1}^1$ for all $n\in \mathbb{Z}$ we get
$\Delta_{11}(\la)\equiv\tilde\Delta_{11}(\la).$ By virtue of (\ref{16}),
this yields
\begin{align}\label{30}
N(\la)\equiv\tilde N(\la).
\end{align}

Define a matrix $P(x,\lambda)$ as
\begin{equation}\label{31}
\begin{array}{l}
P(x,\lambda)(\tilde{v}_1(x,\lambda),\tilde{v}_2(x,\lambda))=(v_1(x,\lambda),v_2(x,\lambda))
\end{array}
\end{equation}
where
$$
P(x,\lambda)=\left(\begin{array}{cc}
p_{11}(x,\lambda)&p_{12}(x,\lambda)\\
p_{21}(x,\lambda)&p_{22}(x,\lambda)\end{array}\right).
$$
Note that
$$
\left(\begin{array}{cc}
\tilde{v}_{11}(x,\lambda)&\tilde{v}_{21}(x,\lambda)\\
\tilde{v}_{12}(x,\lambda)&\tilde{v}_{22}(x,\lambda)\end{array}\right)^{-1}=\left(\begin{array}{cc}
\tilde{v}_{22}(x,\lambda)&-\tilde{v}_{21}(x,\lambda)\\
-\tilde{v}_{12}(x,\lambda)&\tilde{v}_{11}(x,\lambda)\end{array}\right),
$$
thus
\begin{align*}
&p_{11}(x,\lambda)=v_{11}(x,\lambda)\tilde{v}_{22}(x,\lambda)-
\tilde{v}_{12}(x,\lambda)v_{21}(x,\lambda),\\
&p_{12}(x,\lambda)=-v_{11}(x,\lambda)\tilde{v}_{21}(x,\lambda)+\tilde{v}_{11}(x,\lambda)v_{21}(x,\lambda),\\
&p_{21}(x,\lambda)=v_{12}(x,\lambda)\tilde{v}_{22}(x,\lambda)-\tilde{v}_{12}(x,\lambda)v_{22}(x,\lambda),\\
&p_{22}(x,\lambda)=-v_{12}(x,\lambda)\tilde{v}_{21}(x,\lambda)+\tilde{v}_{11}(x,\lambda)v_{22}(x,\lambda).
\end{align*}
Using (\ref{15}) and (\ref{30}), one gets
\begin{align*}
p_{11}(x,\lambda)&=(Z_{11}(x,\la)\tilde Z_{22}(x,\la)-\tilde Z_{12}(x,\la)Z_{21}(x,\la))\\
&+(\tilde N(\la)-N(\la))Z_{11}(x,\la)\tilde Z_{12}(x,\la)\\
&=Z_{11}(x,\la)\tilde Z_{22}(x,\la)-\tilde Z_{12}(x,\la)Z_{21}(x,\la).
\end{align*}
Thus, for each fixed $x\in (0,T),$ the function $p_{11}(x,\lambda)$
is entire in $\la.$ On the other hand, taking (\ref{18}) and (\ref{21}) into
account we calculate for each fixed $x\ge T/2$:
$$
p_{11}(x,\lambda)-1=o(1),\; |\rho|\to\iy,\;\rho\in\Pi_\de.
$$
Also, applying
(\ref{26}) and (\ref{29}), we get
$$
p_{11}(x,\lambda)=O(1),\; |\rho|\to\iy,\;\rho\in G_\de'.
$$
Using the maximum modulus principle and Liouville's theorem for
entire functions, we conclude that
$$
p_{11}(x,\lambda)\equiv 1,\quad x\ge T/2.
$$
Similarly, we obtain
$$
p_{12}(x,\lambda)=p_{21}(x,\lambda)\equiv 0,\ p_{12}(x,\lambda)\equiv 1,\ \quad x\ge T/2.
$$
Together with (\ref{31}) this yields that for $x\ge T/2$,
\begin{align}\label{32}
v_k(x,\la)=\tilde v_k(x,\la),\;
Z_k(x,\la)=\tilde Z_k(x,\la),\; \Omega(x)=\tilde \Omega(x).
\end{align}

Next let us now consider the BVPs $L_1^a$ and $L_{11}^a$ for Eq.(\ref{1}) on
the interval $(0,T)$ with the conditions $U_1^a(y)=V_1(y)=0$ and
$U_1^a(y)=V_2(y)=0,$ respectively. Then, according to Eq.(\ref{10}), the
functions $\Delta_1^a(\la):=-U_1^a(Z_2)$ and $\Delta_{11}^a(\la)
:=U_1^a(Z_1)$ are the characteristic functions of $L_1^a$ and
$L_{11}^a,$ respectively. And
$$
U_1^{a/2}(Z_k)=U_1^{a}(Z_k)-\int_{a/2}^a Z_k(t,\la)d\sigma_1(t),
\quad k=1,2,
$$
hence
\begin{equation}\label{33}
\begin{array}{l}
\Delta_{1}^{a/2}(\la)=
\Delta_{1}^a(\la)+\int_{a/2}^a Z_2(t,\la)d\sigma_1(t),\\
\Delta_{11}^{a/2}(\la)=
\Delta_{11}^a(\la)-\int_{a/2}^a Z_1(t,\la)d\sigma_1(t).
\end{array}
\end{equation}
Let us use (\ref{33}) for $a=T.$ Since $\Delta_{1}^T(\la)=\Delta_{1}(\la),\;
\Delta_{11}^T(\la)=\Delta_{11}(\la),$ it follows from (\ref{32})-(\ref{33}) that
$$
\Delta_{1}^{T/2}(\la)=\tilde\Delta_{1}^{T/2}(\la),\quad
\Delta_{11}^{T/2}(\la)=\tilde\Delta_{11}^{T/2}(\la).
$$
Repeating preceding arguments subsequently for $a=T/2, T/4, T/8,\ldots,$
we conclude that $\Omega(x)=\tilde \Omega(x)$ on $(0,T).$ Theorem 2 is proved.\B

{\textbf{Proof of Theorem 1}}

Define the matrix $R(x,\lambda)$ as
\begin{equation}\label{34}
\begin{array}{l}
R(x,\lambda)(\tilde{\Phi}(x,\lambda),\tilde{\varphi}(x,\lambda))=(\Phi(x,\la),\vv(x,\la))
\end{array}
\end{equation}
where
$$
R(x,\lambda)=\left(\begin{array}{cc}
R_{11}(x,\lambda)&R_{12}(x,\lambda)\\
R_{21}(x,\lambda)&R_{22}(x,\lambda)\end{array}\right).
$$
Note that
$$
\mbox{ det }\left(\begin{array}{cc}
\tilde{\Phi}_{1}(x,\lambda)&\tilde{\varphi}_{1}(x,\lambda)\\
\tilde{\Phi}_{2}(x,\lambda)&\tilde{\varphi}_{2}(x,\lambda)\end{array}\right)=1,
$$
thus
\begin{align*}
&R_{11}(x,\lambda)=\Phi_{1}(x,\lambda)\tilde{\varphi}_{2}(x,\lambda)-
\varphi_{1}(x,\lambda)\tilde\Phi_{2}(x,\lambda),\\
&R_{12}(x,\lambda)=\Phi_{1}(x,\lambda)\tilde{\varphi}_{1}(x,\lambda)-
\varphi_{1}(x,\lambda)\tilde\Phi_{1}(x,\lambda),\\
&R_{21}(x,\lambda)=\Phi_{2}(x,\lambda)\tilde{\varphi}_{2}(x,\lambda)-
\varphi_{2}(x,\lambda)\tilde\Phi_{2}(x,\lambda),\\
&R_{22}(x,\lambda)=\Phi_{2}(x,\lambda)\tilde{\varphi}_{1}(x,\lambda)-
\varphi_{2}(x,\lambda)\tilde\Phi_{1}(x,\lambda).
\end{align*}

Since $\Lambda_1\cap\Xi=\emptyset$ we can infer that $\Lambda_1\cap\Lambda_2=\emptyset$. Otherwise,
if a certain $\lambda\in \Lambda_1\cap\Lambda_2$ then $\lambda\in \Xi$. Thus $\lambda\in \Lambda_1\cap\Xi$;
this leads to a contradiction to the assumption that $\Lambda_1\cap\Xi=\emptyset$. Moreover, Eqs.
$M(\la)=\tilde M(\la)$, $M(\la)=\frac{\Delta_2(\lambda)}{\Delta_1(\lambda)}, \mbox{ and} \ \tilde{M}(\la)=\frac{\tilde{\Delta}_2(\lambda)}{\tilde{\Delta}_1(\lambda)}$
imply that
$$
\Delta_1(\lambda)=\tilde{\Delta}_1(\lambda),\ \Delta_2(\lambda)=\tilde{\Delta}_2(\lambda).
$$
It follows from (\ref{11}) and (\ref{34}) that
\begin{align*}
&R_{11}(x,\la)=\frac{1}{\Delta_1(\la)}
\Big(\psi_1(x,\la)\tilde\vv_2(x,\la)-\tilde\psi_2(x,\la)\vv_1(x,\la)\Big),\\
&R_{12}(x,\la)=\frac{1}{\Delta_1(\la)}
\Big(\psi_1(x,\la)\tilde\vv_1(x,\la)-\tilde\psi_1(x,\la)\vv_1(x,\la)\Big),\\
&R_{21}(x,\la)=\frac{1}{\Delta_1(\la)}
\Big(\psi_2(x,\la)\tilde\vv_2(x,\la)-\tilde\psi_2(x,\la)\vv_2(x,\la)\Big),\\
&R_{22}(x,\la)=\frac{1}{\Delta_1(\la)}
\Big(\psi_2(x,\la)\tilde\vv_1(x,\la)-\tilde\psi_1(x,\la)\vv_2(x,\la)\Big).
\end{align*}
The above equations imply that for each fixed $x,$ the function $R_{11}(x,\lambda)$ is
meromorphic in $\la$ with possible poles only at $\la=\la_{n1}$.
On the other hand, taking (\ref{12}) into account, we also get
\begin{align}\label{35}
R_{11}(x,\lambda)=\frac{1}{\om(\la)}
\Big(\theta_1(x,\la)\tilde\vv_2(x,\la)-\tilde\theta_2(x,\la)\vv_1(x,\la)\Big).
\end{align}
The assumption that $\Lambda_1\cap\Xi=\emptyset$ tells us that the function $R_{11}(x,\lambda)$ is regular at $\la=\la_{n1}$. Thus,
for each fixed $x,$ the function $R_{11}(x,\lambda)$ is entire in $\la.$
Using (\ref{17}) and (\ref{20}), we can obtain for $x\ge T/2:$
$$
R_{11}(x,\lambda)-1=o(1),\quad |\rho|\to\iy,\; \rho\in\Pi_\de.
$$
Also, using (27)-(28), we obtain for $x\geq T/2:$
$$
R_{11}(x,\lambda)=O(1),\quad |\rho|\to\iy,\; \rho\in G_\de.
$$
Therefore, $R_{11}(x,\lambda)\equiv 1$ for $x\geq T/2$. Similarly, we have
$R_{12}(x,\lambda)=R_{21}(x,\lambda)=R_{22}(x,\lambda)-1\equiv 0$ for $x\geq T/2$.
Together with (\ref{14}) and (\ref{34}), it yields
$$
\vv(x,\la)=\tilde\vv(x,\la),\; \psi(x,\la)=\tilde\psi(x,\la),\; \Omega(x)=\tilde{\Omega}(x),\; x\ge T/2.
$$
Also, we obtain
$$
Z_k(x,\la)=\tilde Z_k(x,\la),\quad k=1,2,\quad x\ge T/2.
$$
Since
$$
\vv(x,\la)=U_1(Z_1)Z_2(x,\la)-U_1(Z_2)Z_1(x,\la)
$$
and
$$
\tilde{\vv}(x,\la)=U_1(\tilde{Z}_1)\tilde{Z}_2(x,\la)-U_1(\tilde{Z}_2)\tilde{Z}_1(x,\la)
$$
we have
$$
\vv(x,\la)=\Delta_{11}(\lambda)Z_2(x,\la)+\Delta_{1}(\lambda)Z_1(x,\la)
$$
and
$$
\tilde{\vv}(x,\la)=\tilde{\Delta}_{11}(\lambda)\tilde{Z}_2(x,\la)+\tilde{\Delta}_{1}(\lambda)\tilde{Z}_1(x,\la).
$$
Taking $x=T$ we get
$$
\vv(T,\la)=\Delta_{1}(\lambda),\ \tilde{\vv}(T,\la)=\tilde{\Delta}_{1}(\lambda),\ \vv'(T,\la)=\Delta_{11}(\lambda),\ \tilde{\vv}'(T,\la)=\tilde{\Delta}_{11}(\lambda).
$$
It follows from $\varphi(x,\la)=\tilde\varphi(x,\la)$ for $x\geq T/2$ that
$$
\Delta_{1}(\la)=\tilde\Delta_{1}(\la),\;\Delta_{11}(\la)=\tilde\Delta_{11}(\la).
$$
By Theorem 2 we conclude that $\Omega(x)=\tilde{\Omega}(x)$ on $(0,T).$
Theorem 1 is proved.
\B

\section{Counterexamples}

{\textbf{Example 1} \ (To illustrate that if condition $S$ does not hold then Theorem 1 is false)}

Suppose that $T=\pi,\;U_1(y)=y_1(0),\;U_2(y)=y_1(\pi/2),$ $\Omega(x)=\Omega(x+\pi/2)$ for $x\in(0,\pi/2),$
$\Omega(x)\not\equiv \Omega(\pi-x)$ for $x\in (0,\pi)$.

Take $\tilde \Omega(x):=\Omega(\pi-x)$ for $x\in (0,\pi).$ Then BVP
$\tilde{L}_1$:
Eq.(\ref{1}) with  $\tilde{\Omega}(x)=\Omega(\pi-x)$, $\tilde{y}(x)=y(\pi-x)$, and the conditions $U_1(\tilde{y})=V_1(\tilde{y})=0$;

BVP
$\tilde{L}_2$:
Eq.(\ref{1}) with $\tilde{\Omega}(x)=\Omega(\pi-x)$, $\tilde{y}(x)=y(\frac{3\pi}{2}-x)$, and the conditions $U_2(\tilde{y})=V_1(\tilde{y})=0$;

BVP
$\tilde{L}_0$:
Eq.(\ref{1}) with $\tilde{\Omega}(x)=\Omega(\pi/2-x)$ (also equal to $\Omega(\pi-x)$), $\tilde{y}(x)=y(\frac{\pi}{2}-x)$, and the conditions $U_1(\tilde{y})=U_2(\tilde{y})=0$.

Here $\Delta_1(\lambda)$ is the characteristic function for Eq.(\ref{1}) with
$y_1(0)=0=y_1(\pi)$; $\Delta_2(\lambda)$ is the characteristic function for Eq.(\ref{1}) with
$y_1(\pi/2)=0=y_1(\pi)$;
$\omega(\lambda)$ is the characteristic function for Eq.(\ref{1}) with $y_1(0)=0=y_1(\pi/2)$.
From the above fact the following relations are true:
$$
\Delta_1(\la)=\tilde\Delta_1(\la),\; \Delta_2(\la)=\tilde\Delta_2(\la),\;
\om(\la)=\tilde\om(\la),
$$
and, in view of (\ref{13}), $M(\la)=\tilde M(\la).$

Note that $\Lambda_1$, $\Lambda_2$, and $\Xi$ are sets of zeros for characteristic functions
$\Delta_1(\lambda)$, $\Delta_2(\lambda)$ and $\omega(\lambda)$, respectively. Since $\Omega(x)=\Omega(x+\pi/2)$,
there holds
$\omega(\lambda)=\Delta_2(\lambda)$, i.e. $\Xi=\Lambda_2$. Thus for all $\lambda\in \Xi(=\Lambda_2)$ then it yields
$\lambda\in \Lambda_1$, which implies that $\Xi\cap\Lambda_1\neq \emptyset$.

Now $M(\la)=\tilde M(\la)$ and $\omega(\lambda)=\tilde{\omega}(\lambda)$, but $\Xi\cap\Lambda_1\neq \emptyset$. In Theorem 1
condition $S$ does not hold. In fact, at this case $\Omega(x)\neq \tilde{\Omega}(x):=\Omega(\pi-x)$.
This means, that the specification of $M(\la)$ and $\om(\la)$ does not
uniquely determine the function $\Omega(x)$.

{\textbf{Example 2} \ (To illustrate that even if condition $S$ and $M(\lambda)=\tilde{M}(\lambda)$ hold without
the assumption that $\om(\la)=\tilde\om(\la)$ then Theorem 1 is false)}

Suppose that $T=\pi,\;U_1(y)=y_1(0),\;U_2(y)=y_1(\pi-\al),$ where $\al\in(0,\pi/2).$

Let $\Omega(x)\not\equiv \Omega(\pi-x)$, and $\Omega(x)\equiv 0 $ for $x\in[0,\al_0]\cup
[\pi-\al_0,\pi],$ where $\al_0\in(0,\pi/2).$ If $\al<\al_0$, then $\la_{n2}=
\pi n/\al,\; n\in \mathbb{Z}.$ Choose a sufficiently small $\al<\al_0$ such that
$\Lambda_1\cap\Lambda_2=\emptyset.$ Clearly, such choice is possible. Then
$\Lambda_1\cap\Xi=\emptyset,$ i.e. condition $S$ holds. Otherwise, if a certain $\lambda^*\in \Lambda_1\cap\Xi$, then
$\lambda^*\in \Lambda_1\cap\Lambda_2$; this contradicts to the fact that $\Lambda_1\cap\Lambda_2=\emptyset.$

Take $\tilde \Omega(x):=\Omega(\pi-x)$.
Note that $\Delta_2(\lambda)$ is the characteristic function for Eq.(\ref{1}) with
$$
y_1(\pi-\alpha)=0=y_1(\pi);
$$
and
$\tilde{\Delta}_2(\lambda)$ is the characteristic function for Eq.(\ref{1}) with
$$
y_1(0)=0=y_1(\alpha).
$$
A simple calculation shows $\Lambda_2=\tilde{\Lambda}_2=\left\{\frac{\pi n}{\alpha}, n\in \mathbb{Z}\right\}$.

At this case $\Delta_1(\la)=\tilde\Delta_1(\la),$
$\Delta_2(\la)=\tilde\Delta_2(\la),$ and consequently, $M(\la)=\tilde M(\la).$
Now in Theorem 1 condition $S$ holds, and $M(\la)=\tilde M(\la).$ In fact, in this example, $\Omega(x)\neq \tilde{\Omega}(x):=\Omega(\pi-x)$.
So the true condition $S$ and the specification of $M(\la)$ does not uniquely
determine the functions $\Omega(x)$.

\section{Inverse problem from three spectra}
Fix $a\in(0,T).$ Consider Inverse problem 1 in the case when
$U_1(y):=h_{1,1}y_1(0)+h_{1,2}y_2(0),\; U_2(y):=h_{0,1}y_1(a)+h_{0,2}y_2(a),$ where
$|h_{1,1}|+|h_{1,2}|\neq 0$ and $|h_{0,1}|+|h_{0,2}|\neq 0$. Then the boundary value problems $L_0, L_1, L_2$
take the forms
\begin{align*}
&L'_0:\ \mbox{Eq.}(\ref{1}) \mbox{ with } h_{1,1}y_1(0)+h_{1,2}y_2(0)=h_{0,1}y_1(a)+h_{0,2}y_2(a)=0,\\
&L'_1:\ \mbox{Eq.}(\ref{1}) \mbox{ with } h_{1,1}y_1(0)+h_{1,2}y_2(0)=y_1(T)=0,\\
&L'_2:\ \mbox{Eq.}(\ref{1}) \mbox{ with } h_{0,1}y_1(a)+h_{0,2}y_2(a)=y_1(T)=0.
\end{align*}

Denote by $\Lambda'_j=\{\la'_{nj}\}$ the spectrum of $L'_j\ (j=0,1,2),$ and assume that
$\Lambda'_0\cap\Lambda'_1=\emptyset$ (condition $S'$).

{\bf Theorem 4. }{\it Let condition $S'$ hold. If $\Lambda'_j=\tilde\Lambda'_j$,
$j=0,1,2,$ then $\Omega(x)=\tilde \Omega(x)$ on $(0,T).$}

\smallskip
Comparing Theorem 4 with Theorem 1 we note that $\Lambda_0', \Lambda_1'$ and $\Lambda_2'$ correspond to
$\Xi,\Lambda_1$, and $\Lambda_2$ in Theorem 1, respectively. The theorem is a consequence of Theorem 1.
Theorem 4 for Dirac operator is the analogue to Sturm-Liouville operator with two-point separated boundary conditions,
which was studied by many authors
(see, for example, \cite{17,18}).

\noindent {\bf Acknowledgments.}
The research work of the second author was supported by the Russian Ministry of
Education and Science (Grant 1.1436.2014K), and by Grant
13-01-00134 of Russian Foundation for Basic Research. The first author was supported in
part by the National Natural Science Foundation of China (11171152) and
Natural Science Foundation of Jiangsu Province of China (BK 20141392).

\end{document}